\input amstex
\documentstyle{conm-p}
\NoBlackBoxes
%

\define\BZ{\Bbb Z}

\define\Center{\operatorname{Center}}
\define\tB{\tilde B}
\define\tT{\tilde T}

\define\ep{\endproclaim}
\define \a{\alpha}
\define \be{\beta}

\define \CP{\Bbb C\Bbb P}

\define \CPt{\Bbb C\Bbb P^2}

\define \R{\Bbb R}

\define \Z{\Bbb Z}

\define \BC{\Bbb C}

\define \1{^{-1}}\define \2{^{-2}}

\define \df{\dsize\frac}

\define \edm{\enddemo}

\def\C{\Bbb C}

\def\R{\Bbb R}

\def\Z{\Bbb Z}

\topmatter

\title New invariants for surfaces
\endtitle
\author Mina Teicher
\endauthor
\address Department of Mathematics, Bar-Ilan University, 52900 Ramat-Gan,
Israel\endaddress
 \leftheadtext{M. Teicher}
 \email teicher\@macs.biu.ac.il\endemail

  \abstract We define  a new invariant of surfaces, stable on connected
components of
moduli spaces of surfaces.
The new invariant comes from the polycyclic structure of the fundamental
group of the
complement of a branch curve.
We compute this invariant for a few examples.
Braid monodromy factorizations related to curves is a first step in
computing the fundamental group of the complement of the curve, and thus  we
 indicate the possibility of using braid monodromy factorizations of branch
curves as an invariant of a surface.
  \endabstract
 \thanks The author wishes to thank the Max-Planck Institute for a fruitful and
pleasant stay while writing this paper. This paper was partially supported
by the Emmy Noether Research Institute for  Mathematics, Bar-Ilan
University, and the
Minerva Foundation of Germany.
\endthanks
\endtopmatter
\document

\head{0. Introduction} \endhead

After the remarkable discovery of the Donaldson invariants and the
Sieberg-Witten
invariants \cite{SW}, there was hope that one can use them in order to
distinguish between different connected components of moduli spaces of
surfaces. In our
earlier work we indicated that we believe that these invariants are not
fine enough for
this differentiation and a more geometrical approach is needed. Indeed, in
1997, M.
Manetti \cite{Ma2} produced examples of surfaces which are diffeomorphic
but are not a
deformation of each other; and thus it is clear that a more direct
geometric approach is
needed. We went on to suggest the following distinguishing invariant:

Let $X$ be a complex algebraic surface of general type embedded in
$\CP^N$. Take a generic projection of $X$ to
$\CP^2$ and let
$S_X$ be its branch curve. Clearly
$\pi_1 (\CP^2-S_X)$ is stable on a connected component of moduli spaces of
surfaces.
We believe that these groups can distinguish between different components.
We base our
belief on the structure of such groups which  already have been computed.

\head{1.  History of computations of fundamental groups of complements of
branch
curves}\endhead

Let $X$ be a complex algebraic surface of general type embedded in
$\CP^N$. Let
$f : X \rightarrow \CP^2$ be a generic projection and let
$S_X \in \CP^2$ be its branch curve. Let
$\C^2$ be a big affine piece of
$\CP^2$ s.t.
$S_X$ is transversal to the line at infinity. We denote:
$$\bar G_X = \pi_1(\CP^2-S_X), \quad G_X = \pi_1 (\C^2-S_X).$$

The first computations of
$\bar G_X$ were done by Zariski \cite{Z}. He computed
$\bar G_X$ for $X$ a cubic surface in
$\CP^2$, to be
$\Z_2 \ast \Z_3$ (free product of two finite cyclic groups). The topic was
renewed by
Moishezon in the late 1970's, when he generalized Zariski's result to
$X_n$, a $\deg n$ surface in
$\CP^3$, and proved that in this case
$G_{X_n}$ is the braid group
$B_n$ and
$\bar G_{X_{n}}$ is the braid group over its center,
$B_n/\text{center}$. In fact, Moishezon's result \cite{Mo1} for
$n=3$ coincides with Zariski's result since
$B_3/\text{center} \simeq \Z_2 \ast \Z_3$.
$(B_3 = \langle x,y \bigm| xyx = yxy \rangle, \text{center}(B_3) = \langle
(xy)^3
\rangle$. Thus in
$B_3/\text{center}$,
$x = (xy)^3x = xy \cdot xy \cdot xyx$ and
$y = y(xy)^3 = yxy \cdot (xy)^2 = xyx \cdot (xy)^2.$ Thus
$B_3/\text{center}$ is generated by
$xy$ and
$xyx$ while
$(xy)^3 = 1$ and
$(xyx)^2 = xyx \cdot yxy = (xy)^3 =1$).

The next example was
$V_2$ Veronese of order 2 \cite{MoTe3}. In all the above examples, $G$
contained a free
noncommutative subgroup of two elements (and in other related examples as
in \cite{DOZa}
and \cite{GaTe}). We call such a group ``big''. This gave the basis to the
feeling that
this will always be the case. Since 1991 new examples were computed and
these examples
were not ``big''.

\head{2. The new examples: almost solvable groups} \endhead

To our great surprise the invariants that were computed after 1991 were not
``big''. It
turns out that in all of the new examples,
$\bar G_X = \pi_1(\CP^2-S_X)$ and
$G_X = \pi_1 (\C^2-S_X)$ satisfy the following conditions:
\roster
\item"1)" There exists a quotient of the braid group, namely
$\tilde B_n$, s.t.
$\tilde B_n$ acts on
$G_X$ and
$\bar G_X$.
\item"2)"
$G_X$ and
$\bar G_X$ are not ``big''.
\item"3)"
$G_X$ and
$\bar G_X$ are ``small'': They are almost solvable (or virtually solvable
in another
terminology), i.e., they contain a solvable subgroup of finite index.
\item"4)" Moreover,
$G_X$ and
$\bar G_X$ are extensions of a solvable group by a symmetric group.
\item"5)"
$\bar G_X = G/\text{central element}$.
\endroster

Moreover, in all the new examples we have the following series:
$$1 \triangleleft A_0 \triangleleft A_1 \triangleleft A_2 \triangleleft A_3
\triangleleft G_X \tag 1$$ where
$G/A_3 \simeq S_n$ and
$A_i/A_{i-1}$ is a direct sum of
$\Z$ and one or two finite cyclic groups, to some power
$(\Z\oplus \Z_p \oplus \Z_q)^t$. For example we have the following
computations:

\subhead 2.1. $V_p,$ Veronese of order $p,$\ $ p \ge 3$\endsubhead There
exists a series
$$ 1 \triangleleft A_0 \triangleleft A_1 \triangleleft A_2 \triangleleft A_3
\triangleleft G_{V_p}$$
s.t.
$$\align &G_{V_p}/A_3 \simeq S_{p^2},\  \text{symmetric group  on
$p^2$ elements}\\
\quad\\
&A_3/A_2 \simeq \Z\\
\quad\\
&A_2/A_1 \simeq \left\{ \aligned &(\Z \oplus \Z_3)^{p^2-1} \qquad p \equiv
0 (\text{mod}
\, 3) \\ & \Z^{p^2-1} \qquad\quad\qquad p \not \equiv 0 (\text{mod} \, 3)
\endaligned
\right.\\ \quad\\
&A_1/A_0 \simeq \left\{ \aligned \Z_2 &\qquad p \quad \text{odd} \\ 0
&\qquad p \quad
\text{even} \endaligned \right.\\
\quad\\
&A_0=1\endalign$$
Moreover, we know that
$A_1 $ is the commutant subgroup of
$A_3$, and
$A_1 \subset \Center (G_{V_{p}})$.

There exists a series
$$ 1\triangleleft\bar A_0 \triangleleft \bar A_1 \triangleleft \bar A_2
\triangleleft
\bar A_3
\triangleleft \bar G_{V_{p}} $$
s.t.
$$\align &\bar G_{V_{p}}/\bar A_3 \simeq S_{p^2}\\
\quad\\ &\bar A_3/\bar A_2 \simeq \Z_q,  \quad q = \frac{3p(p-1)} 2\\
\quad\\
&\bar A_2/\bar A_1 \simeq \left\{ \aligned &(\Z \oplus \Z_3)^{p^2-1} \qquad p
\equiv 0 (\text{mod} \, 3) \\ & \Z^{p^2-1} \qquad\ \quad\qquad p \not \equiv 0
(\text{mod}
\, 3)
\endaligned
\right.\\
\quad\\
&\bar A_1/\bar A_0 = \left\{ \aligned \Z_2 &\quad p \quad \text{odd} \\ 0
&\quad p
\quad \text{even} \endaligned \right.\\
\quad\\
&\bar A_0=1\endalign$$

\subhead 2.2. $X_{ab},$ projective embedding of
$\CP^{1} \times \CP^{1}$ w.r.t. to
$| al_1 + bl_2 |$\endsubhead There exists a series
$$1 \triangleleft A_0 \triangleleft A_1 \triangleleft A_2 \triangleleft A_3
\triangleleft
G_{X_{ab}} \quad $$ s.t. for
$n = 2ab = \deg X_{ab}$ we have:
$$\allowdisplaybreaks\align &G_{X_{ab}}/A_3 \simeq S_n\\
\quad\\
&A_3/A_2 \simeq \Z\\
\quad\\
&A_2/A_1 \simeq \left\{ \aligned & ( \Z_2)^{n-1}_0\oplus(\Z_{a-b})^{n-1}
\qquad a,b \quad
\text{even} \\ & \Z_{2(a-b)}^{n-1} \qquad\qquad\qquad\quad\quad
\text{otherwise}
\endaligned
\right.\\
\quad\\
&A_1/A_0 \simeq \left\{ \aligned & \Z_2 \oplus \Z_2 \qquad a,b \ \text{even} \\
& \Z_2 \qquad\qquad\ \text{otherwise} \endaligned \right.\\
\quad\\
& A_0=1\endalign$$ There exists a series
$$1 \triangleleft \bar A_0\triangleleft \bar A_1 \triangleleft \bar A_2
\triangleleft
\bar A_3 \triangleleft \bar G_{X_{ab}}$$ s.t.
$$\aligned &\bar G_{X_{ab}}/\bar A_3 \simeq S_n \\
 &\bar A_3/\bar A_2 \simeq \Bbb
Z_{m_1},
\qquad m_1 = 3ab-a-b \\
 &\bar A_2/\bar A_1 \simeq A_2/A_1 \\
&\bar A_1/\bar A_0 \simeq
A_1/A_0 \\
&\bar A_0=A_0\endaligned$$

\subhead 2.3. Complete intersection of $\deg n$ (not a
hypersurface)\endsubhead There
exists a series

$$1\triangleleft   A_0 \triangleleft A_1 \triangleleft A_2 \triangleleft A_3
\triangleleft G$$ s.t.

$$\aligned &G/A_3 \simeq S_n \\ &A_3/A_2 \simeq \Z \\ &A_2/A_1 \simeq
\Z^{n-1} \\
&A_1/A_0 \simeq \Z_2 \\
&A_0=1\endaligned$$ In this case $G$ is
$\tilde B_n$ itself.

\medskip

There are other computations in progress on Hirzebruch surfaces
(following\linebreak
\cite{MoRoTe} and \cite{FRoTe}), on K3-surfaces (following \cite{CiMT}) and
on toric
varieties.

The proofs are based on our braid monodromy techniques as presented in
detail in
\cite{MoTe4}, \cite{MoTe6}, \cite{MoTe7}, \cite{MoTe8}, \cite{MoTe9},
\cite{MoTe10},
\cite{Te2}, \cite{Ro}, \cite{FRoTe}, and on the Van Kampen Theorem \cite{VK}.

 In fact, all these groups are polycyclic with $A_1\subseteq \Center(G),$\
$A_1=(A_3)'=(A_2)'.$
Group theoretic classification of these groups might answer questions such
as the
following:
\ Does every almost polycyclic $\tilde B_n$-group appear as the fundamental
group of
complements of branch curves?
\ How many non-isomorphic groups of that type appear?
(We transfer the question to a classification problem in group theory.)

The interest in fundamental groups is growing in general, see for example,
\cite{CaTo},
\cite{GaTe},
\cite{Si}, \cite{BoKa}, \cite{L1}, \cite{L2}, \cite{RoTe}, \cite{To}.

\head{3. The mysterious quotient of
$B_n$ that acts on $G_X$ and
$\bar G_X$} \endhead

In this section we bring the definition of the braid group and we
distinguish certain
elements, called {\it half-twists}. Using half-twists, we present Artin's
structure
theorem for the braid group and the natural homomorphism to the symmetric
group. We also
define transversal half-twists and the quotient of
$B_n$ called
$\tilde B_n$, and quote the almost solvability theorem for this group.

\definition{Definition 3.1}  $\underline{\text{Braid group}\ B_n = B_n[D,K]}$

\flushpar  Let $D$ be a closed disc in
$\R^2, K \subset D, K$ finite. Let $B$ be the group of all diffeomorphisms
$\beta$ of $D$ such that
$\beta (K) = K, \beta\bigm|_{\partial D} = \text{Id}\bigm|_{\partial D}$. For
$\beta_1, \beta_2 \in B$, we say that
$\beta_1$ is equivalent to
$\beta_2$ if
$\beta_1$ and
$\beta_2$ induce the same automorphism of
$\pi_1(D-K,u)$. The quotient of $B$ by this equivalence relation is called
{\it the braid
group}
$B_n [D,K]$\ $  (n = \# K)$. The elements of
$B_n[D,K]$ are called {\it braids}.
\enddefinition

\definition{Definition 3.2} $\underline{H(\sigma)\  \text{half-twist
defined by}\
\sigma}$

\flushpar Let
$D,K$ be as above. Let
$a,b \in K, K_{a,b} = K-a-b$ and
$\sigma$ be a simple path in
$D-\partial D$ connected $a$ with $b$ s.t.
$\sigma \cap K = \{a,b\}$. Choose a small regular neighborhood $U$ of
$\sigma$ and an orientation preserving diffeomorphism
$f : \R^2 \rightarrow \C^1$\ $(\C^1$ is taken with the usual ``complex''
orientation)
such that
$f(\sigma) = [-1,1], f(U) = \{z \in \C^1\bigm|   | z | < 2\}$. Let
$\alpha (r), r \ge 0$, be a real smooth monotone function such that
$\alpha (r) = 1$ for
$r \in [0, \frac 3 2]$ and
$\alpha (r) = 0$ for
$r \ge 2$.

  Define a diffeomorphism
$h : \C^1 \rightarrow \C^1$ as follows. For
$z \in \C^1, z = re^{i\varphi}$, let
$h(z) = re^{i(\varphi +\alpha (r)\pi)}$. It is clear that on
$\{z \in \C^1 \bigm| |z | \le \frac 3 2\}, h(z)$ is the positive rotation by
180$^{\circ}$ and that
$h(z) = $ Id\  on
$\{z \in \C^1 \bigm| | z | \ge 2\}$, in particular, on
$\C^1-f(U)$. Considering
$(f \circ h \circ f^{-1})\bigm|_D$ (we always take composition from left to
right), we
get a diffeomorphism of $D$ which interchanges $a$ and $b$ and is the
identity on $D-U$.
Thus it defines an element of
$B_n[D,K]$, called the {\it half-twist defined by}
$\sigma$ and denoted
$H(\sigma)$.\enddefinition
\medskip

  Using half-twists we build a set of generators for
$B_n$.

\definition{Definition 3.3} $\underline{\text{Frame of}\
B_n[D,K]}$

\flushpar Let $D$ be a disc in
$\R^2$. Let
$K = \{a_1,\ldots, a_n\},$ $K \subset D$. Let
$\sigma_1, \ldots, \sigma_{n-1}$ be a system of simple paths in
$D-\partial D$ such that each
$\sigma_i$ connects
$a_i$ with
$a_{i+1}$ and for
$$i,j \in \{ 1,\ldots, n-1\} , i<j, \quad \sigma_i \cap \sigma_j = \cases
\emptyset &\quad \text{if} \quad | i-j | \ge 2 \\ a_{i+1} &\quad \text{if}
\quad
j=i+1. \endcases.$$ Let
$H_i = H(\sigma_i)$. We call the ordered system of half-twists
$(H_1,\ldots, H_{n-1})$ {\it a frame of}
$B_n[D,K]$ {\it defined by}
$(\sigma_1,\ldots, \sigma_{n-1})$, or a frame of
$B_n[D,K]$ for short.
\enddefinition

\demo{Notation 3.4}
$$\aligned &[A,B] = ABA^{-1}B^{-1} \\ &\langle A,B\rangle = ABAB^{-1}
A^{-1} B^{-1} \\
&(A)_B = B^{-1} AB. \endaligned$$
\enddemo

\proclaim{Theorem 3.5. {\rm {(E. Artin's braid group presentation)}}} Let
$\{H_i\}$ be a frame of
$B_n$. Then
$B_n$ is generated by the half-twists
$\{H_i\}$ and all the relations between
$H_1,\ldots, H_{n-1}$ follow from
$$\aligned [H_i,H_j] = 1 \quad &\text{if} \quad | i-j| > 1 , \\
\langle H_i,H_j \rangle = 1 \quad &\text{if} \quad | i-j|= 1, \endaligned$$
$$1 \le i, j \le n-1.$$
\endproclaim

\demo{Proof} \cite{A} (or \cite{MoTe4}, Chapter 5).
\enddemo

\proclaim{Theorem 3.6} Let
$\{H_i\}$ be a frame of
$B_n$. Then
\roster
\item"(i)" for
$n \ge 2$,
$\Center (B_n)\simeq \Bbb Z$ with a generator
$$\Delta^2_n = (H_1 \cdot \ldots \cdot H_{n-1})^n.$$
\item"(ii)"
$B_2 \simeq \Z$ with a generator
$H_i$.
\endroster
\endproclaim

\demo{Proof} \cite{MoTe4}, Corollary V.2.3.
\enddemo

\proclaim{Proposition 3.7} There is a natural homomorphism
$B_n \rightarrow S_n$ (symmetric group on $n$ elements) defined by
$H_i \rightarrow (i,i+1)  $.
\endproclaim

\demo{Proof} Since the transpositions
$\alpha_i = (i,i+1)$, \ $i=1,\dots,n-1,$ satisfy the relations from Artin's
theorem
(3.5), the above is well defined.
\enddemo

\definition{Definition 3.8} $\underline{P_n},$ the pure braid group

\flushpar The kernel of the above homomorphism is
denoted by
$P_n$.
\enddefinition

\remark{Remark 3.9} The transposition
$\alpha_i$ satisfies a relation that
$H_i$ does not satisfy, which is
$\alpha^2_i = 1$. In fact, it is true for any transposition. Under the
above homomorphism
the image of any half-twist is a transposition and thus any square of a
half-twist
belongs to
$\ker (B_n \rightarrow S_n)$ which is
$P_n$.
\endremark

\definition{Definition 3.10} \underbar{Transversal half-twist, adjacent
half-twist,
disjoint half-twist}

\flushpar Let
$\sigma_1$ and
$\sigma_2$ be two paths in $D$ with endpoints in $K$ which do not intersect
$K$ otherwise
(like in 3.2). The half-twists
$H(\sigma_1)$ and
$H(\sigma_2)$ will be called {\it transversal} if
$\sigma_1$ and
$\sigma_2$ intersect transversally in one point which is not an end point
of either of
the
$\sigma_i$'s.

\flushpar The half-twists
$H(\sigma_1)$ and
$H(\sigma_2)$ will be called {\it adjacent} if
$\sigma_1$ and
$\sigma_2$ have one endpoint in common.

\flushpar The half-twists
$H(\sigma_1)$ and
$H(\sigma_2)$ will be called {\it disjoint} if
$\sigma_1$ and
$\sigma_2$ do not intersect.
\enddefinition

\proclaim{Claim 3.11} Disjoint half-twists commute and adjacent half-twists
satisfy the
triple relation
$ABA = BAB$.
\ep

\demo{Proof} By Proposition 3.7 and the fact that every two half-twists are
conjugated to each other.
\enddemo

\definition{Definition 3.12} $\underline{\tilde B_n}$

\flushpar Let
$Q_n$ be the subgroup of
$B_n$ normally generated by
$[X,Y]$ for
$X,Y$ transversal half-twists.
$\tilde B_n$ is the quotient of
$B_n$ modulo
$Q_n$.
\enddefinition

A basic property of
$\tilde B_n$ is the following:

\proclaim{Theorem 3.13}
$\tilde B_n$ is an almost solvable group; there exist a series
$$\aligned &1 \triangleleft \tilde P^{\prime}_n \triangleleft \tilde P_{n,0}
\triangleleft \tilde P_n \triangleleft \tilde B_n \quad \text{s.t.} \\ &\tilde
B_n/\tilde P_n \simeq \Bbb Z^{n-1} \\ &\tilde P_n /\tilde P_{n,0} \simeq \Z
\\ &\tilde
P_{n,0}/
\tilde P^{\prime}_n \simeq \Z^{n-1}\\ &\tilde P^{\prime}_n \simeq \Z_2 \
(\subseteq
\Center (\tilde B_n). \endaligned$$
\endproclaim

\demo{Proof} \cite{Te1}.\edm

\head{4. Applications of Kulikov's proof of Chisini conjecture} \endhead

V. Kulikov proved in 1998 \cite{Ku} that the so-called Chisini conjecture
is true in many
(and in fact for the most important) families of surfaces. Basically, he
proved that if a
curve
$S$ in
$\CP^2$ is a branch curve of a generic projection to
$\CP^2$, then it is the branch curve of only one generic projection. Thus
if we can distinguish different branch curves by fundamental groups of
their complements,
  we will also be able to distinguish between their corresponding surfaces.

\head{5. The new invariant} \endhead

We conjecture that for many classes of embedded (in
$\CP^N)$ surfaces of general type the fundamental group of the complement
of the branch
curve of a generic projection is almost solvable. Moreover, we believe that
like in all
the new examples there exist:
$$\aligned &1 = A_0 \triangleleft \ldots \triangleleft A_{m-1}
\triangleleft A_m
\triangleleft G \\ &1 = \bar A_0 \triangleleft \ldots \triangleleft \bar
A_{m-1}
\triangleleft \bar A_m \triangleleft G\endaligned$$
s.t.
$$\align & G/A_m \simeq \bar G/\bar A_m \simeq S_n \quad n = \deg X\\
 & A_i/A_{i-1} \simeq (\Z_{t_i} \oplus \Z_{s_i} \oplus\Bbb Z^{r_i})^{q_i}\\
& \bar A_i/\bar A_{i-1} \simeq (\Z_{t_{i}'} \oplus \Z_{s_{i}'} \oplus
\Z^{r_i'})^{q^{\prime}_i}\endalign$$ Thus we attach to $X$ the invariants
$$\aligned &n(X) = (t_1,s_1,r_1,q_1, \ldots, t_m, s_m, r_m, q_m, n) \\
&\bar n(X) =
(t^{\prime}_1, s^{\prime}_1, y^{\prime}_1, q^{\prime}_1, \ldots, t^{\prime}_m,
s^{\prime}_m, r^{\prime}_m, q^{\prime}_m, n). \endaligned$$ Clearly, if $X$
and $Y$ are
in the same connected component, these invariants coincide for them and if
they do not
coincide, they are not in the same component.

\example{Examples 5.1}
$$   \alignat 2 n(V_p) &= (2,4,0,1,3,1,1,p^2-1,1,1,1,1,p^2) &&\qquad p
\quad
\text{odd}\\
  &\text{   } &&\qquad p \equiv 0(3) \\
\quad\\
 &= (2,1,0,1,1,1,1,p^2-1, 1,1,1,1,p^2)
&&\qquad p
\quad
\text{odd} \\ &\text{   } &&\qquad p \not \equiv 0(3) \\
\quad\\
\text{   }  &= (1,1,0,0,3,1,1,p^2-1,1,1,1,1,p^2) &&\qquad p \quad
\text{even} \\
&\text{   } &&\qquad p \equiv 0(3) \\
\text{   }\\ &= (1,1,0,0,1,1,1,p^2-1, 1,1,1,1, p^2) &&\qquad p \quad
\text{even} \\
&\text{   } &&\qquad q \not\equiv 0(3) \endalignat$$
$$\alignat2
\bar n(V_p) &= (2,4,0,1,3,1,1,p^2-1, \frac{3p(p-1)}2,
1,0,1,p^2) &&\qquad p
\quad \text{odd} \\  &\text{   } &&\qquad p \equiv 0(3) \\
\quad\\
&= (2,1,0,1,1,1,1,p^2-1,
\frac{3p(p-1)}2, 1,0,1,p^2) &&\qquad p \quad \text{odd} \\ &\text{   }
&&\qquad p \not
\equiv 0(3) \\
\text{   } \\ &= (1,1,0,0,3,1,1,p^2-1,  \frac{3p(p-1)}2, 1,0,1,p^2)
&&\qquad p \quad
\text{even} \\ &\text{   } &&\qquad p \equiv 0(3) \\
\text{   }\\ &= (1,1,0,0,1,1,1,p^2-1,  \frac{3p(p-1)}2, 1,0,1,p^2) &&\qquad
p \quad
\text{even} \\ &\text{   } &&\qquad q \not\equiv 0(3) \\
\quad\\ n(X_{ab}) &= (2,2,0,1,2,a-b,0,2ab-1,
1,1,1,1,2ab) &&\qquad a,b \quad
\text{even} \\ &= (2,1,0,1,2(a-b),1,0,2ab-1,
1,1,1,1, 2ab)  &&\qquad \text{otherwise}  \\
\quad\\
\bar n(X_{ab}) &= (2,2,0,1,2,a-b,0,2ab-1, 3ab-a-b, 1,0,1,2ab)
&&\qquad a,b
\quad \text{even} \\  &=
(2,1,0,1,2(a-b),1,0,2ab-1, 3ab-a-b, 1,0,1, 2ab)  &&\qquad \text{otherwise}
\endalignat$$
$$n(CI)  = (2,1,0,1,1,1,1,n-1,1,1,1,1,n)$$
On the other hand, $\bar n(CI)$  is not completely determined by $n$.
\endexample

\head{6. Example: precise computation of $G$}\endhead

In this section we give the precise statement of $G_X$ for $X$ the Veronese
of order~3.
>From this structure it will be understood how the almost solvability
phenomenon came
about.

Let $G=\pi_1(\Bbb C^2-S)$ for $S$ the branch curve of $V_3\to\Bbb C\Bbb P^2.$

In order to formulate the theorem we need a few definitions.
\bigskip
\definition{Definitions 6.1}
\smallskip
\flushpar$\underline{G_0(9)}$
\smallskip\flushpar
$G_0(9)$ is a $\Bbb Z_2$ extension of a free group on 8 elements.
We take the following model for $G_0(9)\:$

\flushpar Let $G_0(9)$ be generated by $\{g_i\}_{i=1\ i\ne4}^9$ s.t.
$$[g_i,g_j]=\cases
1\quad&T_i, T_j\ \text{are disjoint}\\
\tau\quad & \text{otherwise}\endcases$$
where $\tau^2=1,$\ $\tau\in\Center (G_0(9)).$

\flushpar We take the following action of $\tB_9$ on $G_0(9)$
$$(g_i)_{\tT_k}=\cases g_i\1\tau&\quad k=i\\
g_i&\quad T_i,T_k\ \text{are disjoint}\\
g_ig_k\1&\quad T_i,T_k\ \text{are not orderly adjacent}\\
g_kg_i&\quad \text{otherwise}\endcases$$
\smallskip
 \flushpar  $\underline{\tilde B_9\ltimes G_0(9)}$
\smallskip\flushpar
Consider the semidirect product $\tilde B_9\ltimes G_0(9)$ w.r.t. the
chosen action.
\smallskip
\flushpar$\underline{N_9}$
\smallskip\flushpar
Let $c=[\tT_1^2,\tT_2^2].$

\flushpar Let $\xi_1=(\tT_2\tT_1\tT_2\1)^2\tT_2^{-2}.$

\flushpar Let $N_9\triangleleft\tB_9\ltimes G_0(9)$ be normally generated
by $c\tau\1$
and $(g_1\xi_1\1)^3.$
\smallskip
\flushpar$\underline{G_9}$
\smallskip\flushpar
Let $G_9=\df{\tB_9\ltimes G_0(9)}{N_9}.$
\smallskip
\flushpar  $\underline{\hat\psi_9}$
\smallskip\flushpar Let $\tilde\psi_9$ be the homomorphism $\tB_9\to S_9$
induced from
the standard homomorphism $B_9\to S_9$ (see 3.7).
$\tilde\psi_9$ exists since $[X,Y]\to1$ under the standard homomorphism.
Let $\hat\psi_9: G_9\to S_9$ be defined by the first coordinate
$\hat\psi_9(\a,\be)=\tilde\psi_9(\a).$
\smallskip
\flushpar  $\underline{\psi }$
\smallskip\flushpar
The projection $V_3\to \CPt$, of degree 9, induces a standard monodromy
homomorphism
 $\pi_1(\BC^2-S,*)\to S_9$ which we denote by $\psi .$ \enddefinition

\proclaim{Theorem 6.2}
$G\simeq\tilde B_9\ltimes G_0(9)/N_9$ s.t. $\psi$ is compatible with
$\psi_9.$\ep

\demo{Proof} \cite{MoTe9}, \cite{MoTe10}, \cite{Te2}.\enddemo

$G_9$ is almost polycyclic. More precisely, let\smallskip
\flushpar   $\underline{H_9,H_{9,0},H_9',H_{9,0}'}$
\smallskip\flushpar
Let $Ab: B_9\to\Bbb Z$ be the abelianization of $B_9$  and $B_9$ over its
commutator
subgroup.

\flushpar Let $\widetilde{Ab}:\tB_9\to\BZ$ be a homomorphism induced from
$Ab$ (which
exists since $Ab([X,Y])=1).$

\flushpar Let $\widehat{Ab}: G_9\to\BZ$ be defined by the first coordinate
$\widehat{Ab}(\a,\be)=\widetilde{Ab}(\a).$

\flushpar Let $H_9=\ker\hat\psi_9.$

\flushpar Let $H_{9,0}=\ker \hat\psi_9\cap\ker \widehat{Ab}.$

\flushpar Let $H_9',H_{9,0}'$ be the commutant subgroup of $H_9$ and $H_{9,0},$
respectively.

\proclaim{Proposition 6.3} We have \quad
$1\triangleleft H_{9,0}'\triangleleft H_{9,0}\triangleleft H_9\triangleleft
G_9,$\quad where\quad $G_9/H_9\simeq S_9,$ \linebreak
$H_9/H_{9,0}\simeq\Bbb Z,$\quad
   $H_{9,0}/H'_{9,0}\simeq
(\Bbb Z\oplus\Bbb Z/3\BZ)^8,$\quad $H_{9,0}'=H_9' \simeq \Bbb Z/2 \Bbb Z.$
\endproclaim

\demo{Proof} \ \cite{MoTe10}, Proposition 2.4.\enddemo

In \cite{Te2}, we proved an almost solvability result for the projective
complement.
\head 7. A more basic invariant:  braid
monodromy factorizations related to branch curves
\endhead

The first step in computing the fundamental group of the complement of a
curve is
to compute its braid monodromy. In fact, in order to really realize the
group, one has
to compute braid monodromy factorizations of
$\Delta^2$ (the central element of braid groups related to the curve).

There are many interesting questions which are still open.  They include:

 Is the data in a
braid monodromy factorization of
$\Delta^2$, related to a branch curve,   enough in order to distinguish
between different
connected components of moduli spaces of surfaces? In recent research of V.
Kulikov and
the author \cite{KuTe}, it was proved that if two curves have equivalent braid
monodromy factorizations, the pairs
$(\CP^2,B_i)$ are homeomorphic.
Moreover, if the equivalent braid monodromies are related to branch
curves  then the
associated  surfaces are diffeomorphic.

Will equivalent braid monodromy indicate deformation type within the
algebraic surfaces
category or outside of it?

 Is
$\pi_2$ (the second homotopy group) as a module over
$\pi_1$ needed for this purpose?

How does one  distinguish between two types of braid monodromy
factorizations: those
induced from algebraic curves and those induced from other curves (see
\cite{Mo2})?

How does one determine whether two braid monodromy factorizations are
equivalent?

Can one
use  braid monodromy factorizations to find simplectic invariants (see
\cite{CKTe})?

\medskip

References on diffeomorphism types can be found in \cite{Ma1}, \cite{Ma2}
and \cite{C}.

\enddocument